\pdfoutput=1
\documentclass[11pt,a4paper]{article}
\usepackage[T1]{fontenc}\usepackage{lmodern}
\usepackage[margin=2.5cm]{geometry}\usepackage[hidelinks]{hyperref}
\usepackage{amsmath,amssymb}
\usepackage{graphicx}
\usepackage{booktabs}
\usepackage{array}
\usepackage{orcidlink}
\usepackage{url}
\newtheorem{thm}{Theorem}
\newtheorem{prop}{Proposition}
\newtheorem{lem}{Lemma}
\newtheorem{cor}{Corollary}
\newtheorem{rem}{Remark}
\newenvironment{pf}{\par\noindent\textit{Proof.}\ }{\hfill$\square$\par\medskip}
\newcommand{\dpi}{\Delta\pi}
\newcommand{\vphi}{\varphi}
\title{Analytical modelling of compliance dynamics under fiscal feedback: tipping thresholds, intervention costs and hysteresis}
\author{Luca Meacci$^{1,2}$\,\orcidlink{0000-0001-9341-8593},
Francisco J. Mu\~noz$^{3}$\,\orcidlink{0000-0002-0874-6275}, Juan Carlos Nu\~no$^{3}$\,\orcidlink{0000-0003-4708-980X}\thanks{Corresponding author: \texttt{juancarlos.nuno@upm.es}}\\[1ex]
\small $^1$Dipartimento di Matematica e Informatica ``Ulisse Dini'', Universit\`a degli Studi di Firenze,\\
\small Viale Morgagni 67/a, 50134 Firenze, Italy\\[0.6ex]
\small $^2$Dipartimento di Ingegneria dell'Informazione (DINFO), Universit\`a degli Studi di Firenze,\\
\small Via di S.~Marta 3, 50139 Firenze, Italy\\[0.6ex]
\small $^3$Departamento de Matem\'atica Aplicada, Universidad Polit\'ecnica de Madrid,\\
\small 28040 Madrid, Spain}
\date{}
\begin{document}\maketitle
\begin{abstract}
Compliance with taxes and other collective obligations is often enforced with resources that come from those who comply, so that non-compliance erodes its own deterrent. We develop an exactly solvable dynamical model of compliance under this fiscal feedback. Individuals choose between paying and evading by imitating those who do better, and the budget against evasion is a share of the revenue, split between enforcement and education. The model yields closed-form expressions for the equilibria, the tipping thresholds and their bifurcations, the budgets needed to maintain compliance and to restore it after a collapse, and the recovery time near the threshold. Fiscal feedback alone does not produce tipping points: hysteresis appears only when an instrument that changes individual incentives, such as fines or education, is funded by the tax base. Stronger sanctions reduce both the maintenance and the rescue budget, but beyond a critical level they make the system bistable, so that the choice of sanction is a trade-off between efficiency and resilience, which we quantify with a resilience indicator and a resilience-constrained design problem. The cheapest allocation of the budget between enforcement and education is also obtained in closed form. A sensitivity analysis identifies the parameters that control prevention and recovery costs, and stochastic simulations of individual agents and a slow budget cycle confirm the analytical predictions. The framework provides quantitative guidelines for allocating scarce intervention resources in tax compliance and in other enforcement problems with similar feedbacks.
\end{abstract}

\section{Introduction}

Dodging taxes, claiming undue benefits and riding public transport without a ticket share a structure: individuals enjoy a collective resource while withholding their share of its cost. The response of such behaviour to the likelihood and cost of being caught is studied by the economics of crime and of tax compliance~\cite{Becker1968,AllinghamSandmo1972,Andreoni1998,PickhardtPrinz2014}, and evolutionary game theory has shown how sanctioning institutions can sustain contributions in populations that learn by imitation~\cite{HofbauerSigmund1998,Sigmund2010}. Population-dynamical models of crime~\cite{Nuno2008,Nuno2010,Nuno2011} add the budget: the resources used against cheaters come from the society that cheaters exploit. In the model of Nu\~no, Herrero and Primicerio~\cite{Nuno2010} (hereafter NHP), a share $\vphi$ of the common wealth is spent against cheating, a fraction $s$ of it on education and the rest on police.

This coupling between compliance and enforcement is a feedback loop, and feedback loops are the usual origin of alternative stable states and abrupt transitions~\cite{Scheffer2009}. When both a mostly honest and a mostly dishonest equilibrium are stable, the history of a society matters: after a collapse, restoring the budget that used to sustain compliance is not enough. The question we address is which ingredients are needed for this to happen, and how large the effect is. In models with several social classes and many parameters this question is answered case by case, numerically. Here we take the opposite route and look for a setting small enough that the answer can be written down.

NHP proposed, besides their main three-class model, a two-class model with only payers and cheaters, in which recruitment into cheating is spontaneous; they argued from its isoclines that its steady state is unique and stated that it is asymptotically stable. Like most population models of crime, it prescribes the rates at which individuals move between classes. Here we derive these rates instead from individual decisions, using evolutionary game theory~\cite{HofbauerSigmund1998,Sandholm2010}. Each individual either pays or cheats; the two options have payoffs that depend on the policy and, through the common budget, on what everybody else does; and individuals change option by imitating others who do better~\cite{Schlag1998}. The outcome is again a system of differential equations for the proportions of payers and cheaters, so the game-theoretic description does not replace the population description: it provides its microscopic foundation. What it adds is that every policy instrument enters through its effect on payoffs or on removal, so that one can tell which instruments change behaviour. The resulting model is still two-dimensional and, as we show, exactly solvable. All the results below concern this model and policies that are constant in time.

From a modelling point of view, the problem faced by an agency in charge of compliance is one of allocating scarce resources under feedback. Each year it decides how much to spend and how to divide it between enforcement (audits, inspections, sanctions) and prevention (education, information campaigns), while the resources available depend on how many comply. Tax administration is the main example and the one we use throughout, but the same structure appears whenever enforcement is financed by those who comply: fare inspection in public transport paid from ticket revenue, monitoring of fishing or waste regulations financed by licence and collection fees, or anti-corruption units funded from a public budget that corruption depletes. A useful model should tell the decision maker which instruments can push the system across a threshold, how much it costs to prevent a collapse and to recover from one, how sensitive these costs are to uncertain parameters, and how to allocate a given budget. The model developed here is deliberately simple, so that these answers can be obtained in closed form, verified numerically and summarized in a form that can be used for decisions. Throughout, we call the two behaviours \emph{paying} and \emph{cheating} (evading), and the individuals payers and cheaters.

Our results are the following.
\begin{enumerate}
\item The model reduces exactly to a scalar equation for the cheater fraction, whose zeros and slopes give all equilibria and their stability, and it has no periodic orbits (Section~\ref{sec:reduction}).
\item The cheater-free state is bistable with a cheating state if and only if $\lambda\,|\mathrm{d}\dpi/\mathrm{d}c| > \mu$ at full compliance, where $\lambda$ is the imitation rate, $\dpi$ the payoff advantage of cheating, $c$ the fraction of payers and $\mu$ the spontaneous drop-out rate (Theorem~\ref{thm:criterion}). The removal of cheaters by police financed from the tax base drops out of the criterion exactly, so that in this model police alone, or any budget set as an absolute amount, never produces hysteresis, whereas sanctions and education that act on incentives do (Section~\ref{sec:criterion}).
\item With sanctions, the bistable band, the critical sanction, the cusp and the fold curves are given in closed form, also with spontaneous cheating (Sections~\ref{sec:sanctions} and~\ref{sec:tau0}). Rescuing a collapsed society costs at most a fixed multiple, independent of the sanction, of keeping a compliant one, and the time needed to rescue it with a constant policy diverges as the inverse square root of the excess budget (Section~\ref{sec:rescue}).
\item The cheapest allocation of a maintenance budget between enforcement and education is given in closed form; it includes education only if $\lambda\eta > k(1+\lambda F)$, so that sanctions crowd education out. Stronger sanctions lower both the maintenance and the rescue budget but, above a critical value, create a tipping point; a resilience indicator and a resilience-constrained design problem make this trade-off explicit (Section~\ref{sec:optimal}).
\item Closed-form elasticities identify the parameters that control prevention and recovery costs; individual-based stochastic simulations and a slow budget cycle confirm the analytical predictions; and a decision table summarizes maintenance and rescue budgets and recovery times (Sections~\ref{sec:verification} and~\ref{sec:decision}).
\end{enumerate}

Remark~\ref{rem:threeclass} explains why the conclusions of Section~\ref{sec:criterion} are specific to the two-class structure.

\section{Model}\label{sec:model}

\subsection{The two-class model of NHP}

Consider a population of constant size $N$ with $Y$ cheaters and $N-Y$ payers, and let $W$ be the common wealth. Every payer contributes $a$ per unit time. Ordinary expenses take a share $\theta$ of $W$, and anti-cheating spending a share $\vphi$, split into education ($s\vphi W$) and police, i.e.\ enforcement ($(1-s)\vphi W$). NHP~\cite[appendix]{Nuno2010} proposed
\begin{equation}\label{eq:NHP2}
\dot Y = \frac{\bar\tau_1}{1+\eta_1 s\vphi W}\,(N-Y) - \bigl(\mu + k(1-s)\vphi W\bigr)Y,\qquad
\dot W = a(N-Y) - (\theta+\vphi)W .
\end{equation}
Payers become cheaters at a rate $\bar\tau_1$ reduced by education, and cheaters stop cheating spontaneously (rate $\mu$) or are detected by the police (rate $k(1-s)\vphi W$), after which they pay again. For their parameter values ($\mu=0.01$, $k=1$, $\theta=1$, $\bar\tau_1=0.01$, $\eta_1=100$), NHP argued from the isoclines of (\ref{eq:NHP2}) that its equilibrium is unique and stated that it is asymptotically stable. Proposition~\ref{prop:reduction} below makes the stability statement rigorous for the more general model introduced next.

\subsection{A simple evolutionary game}\label{sec:whygame}

Model (\ref{eq:NHP2}) is a population model in the usual sense: it states, at the level of the population, how fast individuals move from one class to another. Its rates are phenomenological; they summarize behaviour without describing why individuals cheat, and a policy can only be represented by its assumed effect on a rate. With the aim of formulating a model in which these rates follow from the decision to cheat itself, we describe the population as playing a simple evolutionary game.

The decision to cheat is naturally a strategic one. Whether cheating pays for an individual depends on what the others do, because detection and education are financed by the taxes of those who pay: when more people cheat, the budget shrinks, the chance of being caught falls and education loses strength, so that cheating becomes more attractive for everybody. This interdependence between individual payoffs, which is the fiscal feedback described in the introduction, is what makes the situation a game.

A classical analysis of this game would assume that every taxpayer knows the probability of detection, the sanction and the behaviour of the others, and chooses a best response; its outcome would be a Nash equilibrium. Such a description says which states are consistent with rational choices, but not how a society moves from one state to another, which is precisely what tipping points and hysteresis are about. It also asks a great deal of individuals who, in practice, know the probability of an audit only vaguely and learn about it largely from the experience of people around them, as the empirical and agent-based literature on tax compliance repeatedly finds~\cite{Andreoni1998,PickhardtPrinz2014}.

Evolutionary game theory~\cite{HofbauerSigmund1998,Sandholm2010} keeps the strategic structure, with strategies whose payoffs depend on the state of the population, but replaces optimization by a simple rule for changing strategy, here the imitation of individuals who are doing better. Successful behaviour spreads and unsuccessful behaviour declines, and the result is a system of differential equations for the proportions of individuals using each strategy, which describes the transient dynamics as well as the rest points. The two viewpoints remain connected, since under imitation every Nash equilibrium is a rest point and every stable rest point is a Nash equilibrium~\cite{HofbauerSigmund1998}.

In this sense evolutionary game theory does not replace a description by differential equations; it derives one from individual behaviour. In our setting each individual either pays ($P$) or cheats ($C$), the payoffs of the two options depend on the policy and on the common wealth (Section~\ref{sec:payoffs}), and individuals revise their choice by imitation (Section~\ref{sec:derivation}). In return, the rates of the population model acquire a behavioural meaning and each policy instrument enters in a definite way: fines and education change payoffs, and hence the decisions of imitators, while detection also removes cheaters directly. The distinction between these two channels is the key to the results of Section~\ref{sec:criterion}.

\subsection{Payoffs}\label{sec:payoffs}

We write $y=Y/N$ for the fraction of cheaters, $c=1-y$ for the fraction of payers, and $w=W/(aN)$, so that $w=1$ is the yearly contribution of a society in which everybody pays. Payoffs are measured per year in units of the yearly tax of a payer. Let $B(w)$ be the benefit that an individual draws from what the common wealth provides. A payer receives this benefit and pays the tax:
\[
\pi_P=B(w)-1 .
\]
A cheater receives the same benefit and saves the tax, but bears two costs. Detection is a random event with rate $d(w)=k(1-s)\vphi w$ per year (the police budget is $(1-s)\vphi$ times the wealth, and $k$ is its efficiency); each detection costs the cheater a sanction $F\ge0$, so the expected sanction per year is $Fd(w)$, as in the economic theory of deterrence~\cite{Becker1968}. In addition, cheating carries a moral or psychological cost $m(w)\in[0,1)$, which education raises; following NHP, who let education reduce recruitment by the factor $1/(1+\eta_1s\vphi W)$, we take $1-m(w)=g(w)=1/(1+\eta s\vphi w)$, with $\eta\ge0$. Thus
\[
\pi_C=B(w)-m(w)-Fd(w),
\]
and cheating is more profitable than paying by
\begin{equation}\label{eq:payoff}
\dpi(w)=\pi_C-\pi_P = g(w) - F d(w),\qquad g(w)=\frac{1}{1+\eta s\vphi w},\qquad d(w)=k(1-s)\vphi w .
\end{equation}
The benefit $B(w)$ cancels: everybody enjoys the public good whether paying or not, which is the social dilemma. We have also neglected the effect of a single individual's tax on $w$, which is of order $1/N$. Without education and sanctions, $\dpi=1$ and cheating always pays.

\subsection{From individual revisions to the population equations}\label{sec:derivation}

We now derive the equation for the fraction of cheaters from the behaviour of individuals, in four steps. Consider a population of $N$ individuals, $Y$ of whom cheat.

\emph{Step 1: the individual rule.} Each individual reconsiders its strategy from time to time, at the times of a Poisson process of rate one per year, independently of the others. At a revision, the individual observes another member of the population, chosen at random, together with that member's payoff. If the other uses the same strategy, nothing changes. If the other uses the other strategy, the reviser adopts it with probability
\begin{equation}\label{eq:schlag}
p_{\rm switch}=\lambda\,\bigl[\pi_{\rm other}-\pi_{\rm own}\bigr]_+,\qquad [x]_+=\max(x,0),
\end{equation}
where $\lambda>0$ is small enough for $p_{\rm switch}\le1$ ($\lambda\le1/\max|\dpi|$). Thus an individual never copies a strategy that does worse, and copies a better one more often the larger the advantage. This is the \emph{proportional imitation rule} of Schlag~\cite{Schlag1998}, who showed that, among rules that use only one's own payoff and that of one sampled individual, rules of this proportional form are those that increase the expected payoff of the population in every situation, and that in a large population they produce the replicator dynamics. Since $\pi_C-\pi_P=\dpi$, a payer who meets a cheater switches with probability $\lambda[\dpi]_+$, and a cheater who meets a payer with probability $\lambda[-\dpi]_+$.

\emph{Step 2: other events.} Three further events change strategies without comparing payoffs: a payer may start cheating spontaneously, at a rate $\tau_0g(w)$ that education reduces as in (\ref{eq:NHP2}) (``exploration'', or mutation in biological language); a cheater may stop spontaneously, at rate $\mu$; and a cheater may be detected, at rate $d(w)$, after which the sanction is paid and the individual resumes paying taxes.

\emph{Step 3: expected change in a short time.} In a short interval $(t,t+\Delta t)$, each of the $N-Y$ payers revises with probability $\Delta t$, meets a cheater with probability $Y/N$ and then switches with probability $\lambda[\dpi]_+$; it also starts cheating spontaneously with probability $\tau_0g\,\Delta t$. Each of the $Y$ cheaters revises with probability $\Delta t$, meets a payer with probability $(N-Y)/N$ and switches with probability $\lambda[-\dpi]_+$, or stops spontaneously or is detected with probabilities $\mu\Delta t$ and $d\,\Delta t$. Up to terms of order $\Delta t^2$, the expected change $\Delta Y$ in the number of cheaters satisfies
\begin{equation}\label{eq:expected}
\frac{\mathbb E[\Delta Y]}{\Delta t}
=\lambda\frac{(N-Y)\,Y}{N}\Bigl([\dpi]_+-[-\dpi]_+\Bigr)+(N-Y)\,\tau_0g-(\mu+d)\,Y ,
\end{equation}
where the first term collects payers who copy a cheater and cheaters who copy a payer, the second spontaneous entry, and the last drop-out and detection.
Dividing by $N$ and writing $y=Y/N$, $c=1-y$, and since $[x]_+-[-x]_+=x$, the two imitation flows combine into a single net flow $\lambda yc\,\dpi$: payers become cheaters when cheating pays more ($\dpi>0$), and cheaters become payers when it pays less ($\dpi<0$).

\emph{Step 4: large population.} The process $Y(t)$ is a Markov jump process whose jumps are of size $\pm1$. As $N\to\infty$, the fraction $y(t)=Y(t)/N$ converges, on every finite time interval, to the solution of the differential equation whose right-hand side is the expected change per unit time computed in Step~3, with fluctuations of order $N^{-1/2}$ (Kurtz's theorem~\cite{Kurtz1970}). The same holds for the wealth, which obeys the balance $\dot W=a(N-Y)-(\theta+\vphi)W$ of (\ref{eq:NHP2}) and, after division by $aN$, $\dot w=c-(\theta+\vphi)w$.

The result is
\begin{equation}\label{eq:model}
\begin{aligned}
\dot y &= \underbrace{\tau_0\, g(w)\, c}_{\text{spontaneous entry}} + \underbrace{\lambda\, y\, c\, \dpi(w)}_{\text{imitation}} - \underbrace{\mu y}_{\text{drop-out}} - \underbrace{d(w)\, y}_{\text{detection}},\\
\dot w &= c - \rho w,\qquad \rho := \theta+\vphi,\qquad c=1-y .
\end{aligned}
\end{equation}
The imitation term is the replicator equation of the two-strategy game: with mean payoff $\bar\pi=y\pi_C+c\pi_P$, one has $y(\pi_C-\bar\pi)=yc\,\dpi$, so imitation makes a strategy grow in proportion to its frequency and to its payoff advantage over the average~\cite{TaylorJonker1978,HofbauerSigmund1998}. The factor $yc$ is the probability that a revision brings together individuals with different strategies; it vanishes when everybody uses the same strategy, since then there is nobody different to imitate. The policy is the triple $(\vphi,s,F)$, held constant in this paper. Fines collected are not added to $W$. The state space $\{0\le y\le1,\ w\ge0\}$ is forward invariant, and $w(t)\le\max\{w(0),1/\rho\}$.

\begin{prop}\label{prop:equiv}
For $\lambda=0$, system (\ref{eq:model}) is the two-class model (\ref{eq:NHP2}) of NHP in normalized variables, with $\tau_0=\bar\tau_1$ and $\eta$, $k$ rescaled by $aN$.
\end{prop}

\begin{pf}
Divide the first equation of (\ref{eq:NHP2}) by $N$ and the second by $aN$, and write $\eta_1 W = \eta_1 aN\,w$ and $kW=kaN\,w$.
\end{pf}

Thus the NHP two-class model is the limit in which cheating spreads only spontaneously. With $\lambda>0$, deterrence has two distinct effects in (\ref{eq:model}): the term $-dy$ removes cheaters, and the term $-\lambda yc\,Fd$ makes cheating less attractive to imitators. Section~\ref{sec:criterion} shows that only effects of the second kind can create hysteresis in this model.

\paragraph{Parameters.} Table~\ref{tab:params} lists the variables and parameters of the model. Unless stated otherwise we use the values of NHP: $\mu=0.01$, $k=1$, $\theta=1$ (time in years), and the imitation rate $\lambda=0.1$, equal to the contagion rate of the three-class NHP model. Spontaneous entry $\tau_0$ is either zero, which isolates the imitation channel, or small ($10^{-3}$, or $10^{-2}$, the value $\bar\tau_1$ of NHP's appendix). For education we use $\eta=100$, the value of NHP's appendix.

\begin{table}[t]
\centering
\caption{Variables and parameters of model (\ref{eq:model}). Rates are per year; budget shares are fractions of the common wealth.}
\label{tab:params}
\small
\begin{tabular}{@{}l >{\raggedright\arraybackslash}p{0.56\textwidth} >{\raggedright\arraybackslash}p{0.2\textwidth}@{}}
\toprule
Symbol & Meaning & Reference value\\
\midrule
$y$, $c=1-y$ & fractions of cheaters and payers & --\\
$w$ & common wealth, in units of the yearly revenue at full compliance & --\\
$\vphi$ & budget share spent against cheating & decision variable\\
$s$ & fraction of that budget spent on education & decision variable\\
$F$ & sanction per detection, in years of tax & decision variable\\
$\theta$ & budget share of ordinary expenses & $1$\\
$k$ & efficiency of enforcement & $1$\\
$\eta$ & efficiency of education & $100$\\
$\lambda$ & imitation rate & $0.1$\\
$\mu$ & spontaneous drop-out rate of cheaters & $0.01$\\
$\tau_0$ & spontaneous entry rate into cheating & $0$, $10^{-3}$, $10^{-2}$\\
$D$, $E$ & enforcement and education strength at full compliance, Eq.~(\ref{eq:DE}) & derived\\
\bottomrule
\end{tabular}
\end{table}

\section{Reduction and global dynamics}\label{sec:reduction}

At an equilibrium of (\ref{eq:model}), $w=c/\rho$. Substituting into the first equation defines the scalar function
\begin{equation}\label{eq:G}
G(y) := \dot y\big|_{w=(1-y)/\rho}
= \tau_0\,\hat g(c)\,c + y\,h(c),\qquad
h(c) := \lambda c\,\widehat{\dpi}(c) - \mu - D c ,
\end{equation}
where hats denote functions evaluated at $w=c/\rho$:
\begin{equation}\label{eq:DE}
\hat g(c)=\frac{1}{1+Ec},\qquad \widehat{\dpi}(c)=\hat g(c)-FDc,\qquad
D:=\frac{k(1-s)\vphi}{\theta+\vphi},\qquad E:=\frac{\eta s\vphi}{\theta+\vphi}.
\end{equation}
The parameter $D$ is the detection rate per individual in a fully compliant society, and $E$ the corresponding strength of education. At equilibrium the detection rate is $d=Dc$: it is proportional to compliance because the budget is a share of a wealth produced by payers. The function $h(c)$ is the per-capita growth rate of cheaters when $\tau_0=0$.

\begin{prop}\label{prop:reduction}
The equilibria of (\ref{eq:model}) are the points $(y^*,(1-y^*)/\rho)$ with $G(y^*)=0$. An equilibrium is asymptotically stable if $G'(y^*)<0$ and a saddle if $G'(y^*)>0$. Moreover, (\ref{eq:model}) has no periodic orbits in $0<y<1$, $w>0$, and if its equilibria are hyperbolic every trajectory converges to an equilibrium.
\end{prop}

\begin{pf}
The first statement follows from $\dot w=0\iff w=c/\rho$. Write $f_1,f_2$ for the right-hand sides of (\ref{eq:model}), and $J_{11}=\partial_y f_1$, $J_{12}=\partial_w f_1$ at the equilibrium. Since $G(y)=f_1(y,(1-y)/\rho)$, $G'=J_{11}-J_{12}/\rho$. The Jacobian is $J=\bigl(\begin{smallmatrix}J_{11}&J_{12}\\-1&-\rho\end{smallmatrix}\bigr)$, with $\det J=-\rho J_{11}+J_{12}=-\rho\,G'$ and $\operatorname{tr}J=G'+J_{12}/\rho-\rho$. Because $g$ decreases and $d$ increases with $w$,
\[
J_{12}=\tau_0 c\,g'(w)+\lambda y c\,\bigl(g'(w)-Fd'(w)\bigr)-d'(w)\,y\le0 .
\]
Hence $G'<0$ implies $\det J>0$ and $\operatorname{tr}J<0$, and $G'>0$ implies $\det J<0$.
For the absence of periodic orbits, take the Dulac function $\beta=1/(yc)$ on $\Omega=\{0<y<1,\,w>0\}$. Then
\[
\frac{f_1}{yc}=\frac{\tau_0 g}{y}+\lambda\dpi-\frac{\mu+d}{c},\qquad \frac{f_2}{yc}=\frac{c-\rho w}{yc},
\]
so that
\[
\partial_y\Bigl(\frac{f_1}{yc}\Bigr)+\partial_w\Bigl(\frac{f_2}{yc}\Bigr)
=-\frac{\tau_0 g}{y^2}-\frac{\mu+d}{c^2}-\frac{\rho}{yc}<0 .
\]
By the Bendixson--Dulac criterion~\cite{Perko2001} there are no periodic orbits in the simply connected domain $\Omega$, nor cycles of orbits connecting saddles. Trajectories are bounded, so by the Poincar\'e--Bendixson theorem every $\omega$-limit set consists of equilibria and, the equilibria being finitely many and hyperbolic, of a single equilibrium.
\end{pf}

Proposition~\ref{prop:reduction} reduces the bifurcation analysis to the scalar function $G$. The wealth relaxes at rate $\rho\approx1$ per year, much faster than the cheater fraction (rates of order $\lambda$ and $\mu$), so $\dot y\approx G(y)$ is also a good approximation of the transient dynamics; we use this in Section~\ref{sec:rescue}.

\section{When does the fiscal feedback create hysteresis?}\label{sec:criterion}

We first set $\tau_0=0$, so that the cheater-free state $y=0$ is an equilibrium and $G(y)=y\,h(1-y)$. Interior equilibria are the zeros of $h$ in $(0,1)$; since $G'(y)=h(c)-y\,h'(c)$, an interior equilibrium is stable when $h'(c)>0$ and a saddle when $h'(c)<0$, and the cheater-free state is stable when $h(1)<0$.

To expose the structure, we state the criterion for general functions. Suppose that, at equilibrium, both $d$ and $\dpi$ depend on compliance through smooth functions $d(c)\ge0$ and $\dpi(c)$, and let $h(c)=\lambda c\,\dpi(c)-\mu-d(c)$.

\begin{thm}\label{thm:criterion}
Let $\tau_0=0$ and let a policy parameter cross a value at which $h(1)=0$ transversally. At this value the cheater-free state undergoes a transcritical bifurcation. If
\begin{equation}\label{eq:criterion}
\sigma := \mu + \lambda\,\dpi'(1) - \bigl(d'(1)-d(1)\bigr) < 0 ,
\end{equation}
the bifurcation is subcritical: on the side where the cheater-free state is stable there is a saddle with $c<1$ and a stable cheating equilibrium, so the society is bistable. If $\sigma>0$ it is supercritical: the cheater-free state is the only attractor where it is stable, and a branch of stable cheating states with $c<1$ emerges where it is unstable.
\end{thm}

\begin{pf}
Differentiating, $h'(1)=\lambda\dpi(1)+\lambda\dpi'(1)-d'(1)$. At the bifurcation $h(1)=0$, that is, $\lambda\dpi(1)=\mu+d(1)$, hence $h'(1)=\sigma$. Near $c=1$, $h(c)\approx h(1)+\sigma(c-1)$, so the interior equilibrium bifurcating from $y=0$ is at $c_r=1-h(1)/\sigma$. If $\sigma<0$ and $h(1)<0$ (cheater-free state stable), then $c_r<1$ and $h'(c_r)<0$: it is a saddle. Moreover $h(0)=-\mu-d(0)<0$ and $h>0$ just below $c_r$, so $h$ has another zero in $(0,c_r)$ at which it increases: a stable cheating equilibrium. If $\sigma>0$, $c_r<1$ requires $h(1)>0$, and the bifurcating equilibrium has $h'(c_r)>0$, so it is stable.
\end{pf}

The criterion can be read as follows. Near full compliance, a fall in compliance weakens deterrence through two channels: incentives, at a rate $\lambda|\dpi'(1)|$ (cheating becomes more profitable), and removal, when the detection rate falls faster than compliance itself, which contributes $d'(1)-d(1)$. Hysteresis appears when the sum exceeds the rate $\mu$ at which cheaters quit on their own.

In the model (\ref{eq:model}), $d(c)=Dc$, so $d'(1)-d(1)=0$: \emph{the mechanical effect of police financed by the tax base cancels exactly}. The reason is that police removal, $Dcy$, and imitation, $\lambda y c\,\dpi$, have the same form; without an effect on incentives, police only lowers the effective imitation rate from $\lambda$ to $\lambda-D$. Writing $\widehat{\dpi}$ from (\ref{eq:DE}), the criterion becomes
\begin{equation}\label{eq:criterion2}
\lambda\Bigl(FD+\frac{E}{(1+E)^2}\Bigr) > \mu \qquad\text{at the threshold}\quad \frac{\lambda}{1+E}=\mu+D(1+\lambda F).
\end{equation}

\begin{lem}\label{lem:concave}
For the model (\ref{eq:model}), $h$ is concave on $[0,1]$. Hence for $\tau_0=0$ there are at most two interior equilibria and at most two attractors.
\end{lem}

\begin{pf}
$h(c)=\lambda c/(1+Ec)-\lambda FDc^2-Dc-\mu$, and each term is concave in $c\ge0$.
\end{pf}

Four consequences follow.

\begin{cor}[police alone]\label{cor:police}
If $F=0$ and $E=0$, then $h(c)=(\lambda-D)c-\mu$ is linear, and for every policy there is a unique attractor: $y^*=0$ if $D\ge\lambda-\mu$, and $y^*=1-\mu/(\lambda-D)$ otherwise. For $\tau_0>0$, $G$ is a quadratic with exactly one zero in $(0,1)$, which is globally attracting.
\end{cor}

\begin{pf}
The first part is immediate. For $\tau_0>0$, $G(0)=\tau_0>0$ and $G(1)=-\mu<0$, and a quadratic with a sign change on $(0,1)$ has exactly one zero there; it attracts every trajectory by Proposition~\ref{prop:reduction}.
\end{pf}

Corollary~\ref{cor:police} extends the uniqueness found by NHP for their two-class model ($\lambda=0$) to any imitation rate: in the two-class setting, imitation alone does not create tipping points, and neither does police financed by the tax base.

\begin{cor}[budget set as an absolute amount]\label{cor:money}
If the resources spent against cheating are a fixed amount $M$ instead of a share of wealth, so that $d=k(1-s)M$ and $g=1/(1+\eta sM)$ do not depend on $c$, then $h$ is linear in $c$ and the attractor is unique for every policy and every $F$.
\end{cor}

In the two-class model, hysteresis therefore requires both a fiscal feedback and a deterrent that acts on incentives. The next two corollaries give the two instances of the latter.

\begin{cor}[sanctions]\label{cor:sanction}
Let $E=0$. The cheater-free state is bistable with a cheating state for some police budget if and only if
\begin{equation}\label{eq:Fc}
F>F_c:=\frac{\mu}{\lambda(\lambda-2\mu)},\qquad\text{which requires}\quad \lambda>2\mu .
\end{equation}
\end{cor}

\begin{pf}
At the threshold $D_T=(\lambda-\mu)/(1+\lambda F)$, (\ref{eq:criterion2}) reads $\lambda F(\lambda-\mu)>\mu(1+\lambda F)$, i.e.\ $\lambda F(\lambda-2\mu)>\mu$. The global statement follows from the explicit band in Proposition~\ref{prop:band}.
\end{pf}

\begin{cor}[education]\label{cor:education}
Let $F=0$. The bifurcation of the cheater-free state is subcritical if and only if $\lambda E/(1+E)^2>\mu$ at the threshold. Since $E/(1+E)^2\le1/4$, education alone can create hysteresis only if $\lambda>4\mu$; for $\lambda=0.1$, $\mu=0.01$ this happens for $E\in(0.127,\,7.87)$.
\end{cor}

With $\eta=100$ and $s=0.5$, for example, $E\approx2.2$ at the threshold, and the society is bistable for $\vphi\in[0.0452,0.0526]$. The mechanism is the same as for sanctions: the moral cost of cheating that education creates is financed by the tax base, so it fades when compliance falls.

\begin{rem}\label{rem:threeclass}
The cancellation behind Corollaries~\ref{cor:police} and~\ref{cor:money} rests on two features of the two-class model: wealth at equilibrium is proportional to the number of payers, so that the detection rate has elasticity exactly one with respect to compliance, and every payer can be recruited into cheating. Neither holds in general. In the three-class NHP model, rich payers contribute more than poor ones and the class composition itself depends on the education budget, so wealth, and with it detection, need not be proportional to compliance, and the pool of payers who can become cheaters changes with wealth. In Theorem~\ref{thm:criterion} the first effect enters through $d'(1)-d(1)$ and the second through an additional factor in the imitation term; either can make $\sigma$ negative without any effect of policy on incentives. The statements of this section about police and about absolute budgets are therefore specific to the two-class structure.
\end{rem}

\section{Sanctions: the bistable band in closed form}\label{sec:sanctions}

Let $\tau_0=0$ and $E=0$ (for instance $s=0$, or $\eta=0$). Then
\begin{equation}\label{eq:hquad}
h(c)=-\lambda FD\,c^2+(\lambda-D)\,c-\mu ,
\end{equation}
and everything can be computed explicitly.

\begin{prop}\label{prop:band}
Let $\tau_0=0$, $E=0$ and $F>F_c$. The society is bistable if and only if $D_T<D<D_{SN}$, where
\begin{equation}\label{eq:band}
D_T=\frac{\lambda-\mu}{1+\lambda F},\qquad
D_{SN}=\lambda\Bigl(\sqrt{1+F\mu}-\sqrt{F\mu}\Bigr)^2 .
\end{equation}
For $D>D_{SN}$ the cheater-free state is the only attractor, and for $D<D_T$ it is unstable and the cheating state is the only attractor. The stable cheating state and the saddle have compliance
\begin{equation}\label{eq:branches}
c_\mp=\frac{\lambda-D\mp\sqrt{(\lambda-D)^2-4\lambda FD\mu}}{2\lambda FD},
\end{equation}
and at $D=D_{SN}$ they merge at $c_v=2\mu/(\lambda-D_{SN})$. At $F=F_c$ the band shrinks to the point $D_T=D_{SN}=\lambda-2\mu$, where the transcritical and the saddle-node bifurcations meet. In terms of the budget share, $\vphi=\theta D/(k(1-s)-D)$.
\end{prop}

\begin{pf}
The zeros of (\ref{eq:hquad}) are (\ref{eq:branches}); $c_-$ is stable and $c_+$ a saddle, since $h'(c_-)>0>h'(c_+)$. They exist when the discriminant is positive, i.e.\ for $D<D_{SN}$, the smaller root of $(\lambda-D)^2=4\lambda F\mu D$, which is (\ref{eq:band}). The cheater-free state is stable if $h(1)<0$, i.e.\ $D>D_T$. Bistability requires in addition $0<c_-<c_+<1$. Given $h(1)<0$ and $h(0)<0$, this holds if and only if the vertex $(\lambda-D)/(2\lambda FD)$ of the parabola lies in $(0,1)$. The vertex decreases with $D$, and at $D=D_T$ it lies below $1$ if and only if $\lambda-D_T<2\lambda FD_T$, which is the condition $F>F_c$ of Corollary~\ref{cor:sanction}. At the fold the vertex is $c_v=2\mu/(\lambda-D_{SN})$, and a direct computation shows that $D_T=D_{SN}=\lambda-2\mu$ (so $c_v=1$) exactly when $F=F_c$, and that $D_T<D_{SN}<\lambda-2\mu$ for $F>F_c$. For $F<F_c$ the double root lies at $c_v>1$, outside the state space, and the bifurcation of the cheater-free state is supercritical.
\end{pf}

\begin{figure}[t]
\centering
\includegraphics[width=0.7\textwidth]{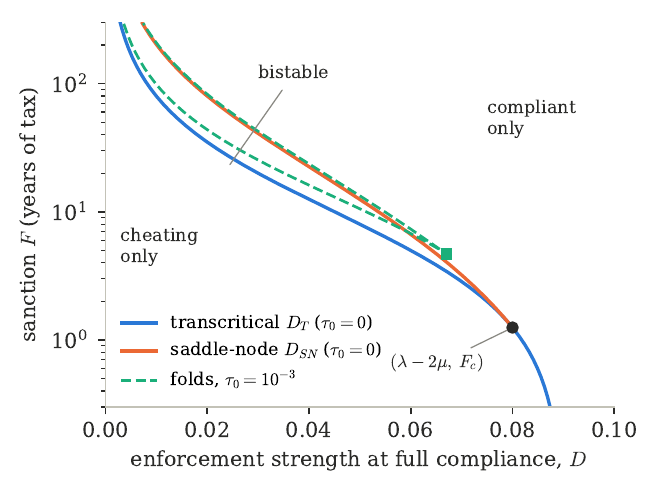}
\caption{Bifurcation set in the plane of enforcement strength $D$ and sanction $F$ ($E=0$). Solid: transcritical line $D_T$ and saddle-node line $D_{SN}$ from (\ref{eq:band}), for $\tau_0=0$; they meet at $(\lambda-2\mu,F_c)$. Dashed: the two fold lines for $\tau_0=10^{-3}$, from (\ref{eq:foldparam}), meeting at the cusp (\ref{eq:cusp}) (square). The budget share is $\vphi=\theta D/(k(1-s)-D)$. Parameters $\lambda=0.1$, $\mu=0.01$.}
\label{fig:phase}
\end{figure}

\begin{figure}[t]
\centering
\includegraphics[width=\textwidth]{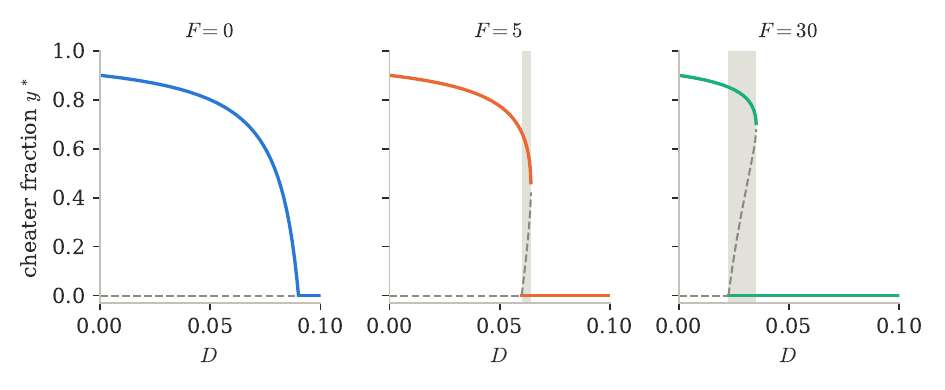}
\caption{Equilibrium cheater fraction versus enforcement strength $D$ for $F=0$, $5$ and $30$ ($\tau_0=0$, $E=0$). Solid: stable; dashed: unstable. Shaded: bistable band (\ref{eq:band}). Without sanction the bifurcation of the cheater-free state is supercritical; above $F_c=1.25$ it is subcritical.}
\label{fig:bif}
\end{figure}

Figure~\ref{fig:phase} shows the band in the $(D,F)$ plane and Fig.~\ref{fig:bif} the equilibria as functions of $D$ for three sanctions; Table~\ref{tab:band} gives numerical values, which agree with a direct continuation of the two-dimensional system (\ref{eq:model}) to the precision of the continuation grid.

\begin{table}[t]
\centering
\caption{Bistable band for $\tau_0=0$, $E=0$, $s=0$, from (\ref{eq:band}). $y_-$: cheater fraction of the stable cheating state at the lower edge ($D=D_T$) and at the fold ($D=D_{SN}$). $R=D_{SN}/D_T$.}
\label{tab:band}
\begin{tabular}{rcccccccc}
\toprule
$F$ & $D_T$ & $D_{SN}$ & $\vphi_T$ & $\vphi_{SN}$ & $R$ & $y_-(D_T)$ & $y_-(D_{SN})$\\
\midrule
2   & 0.0750 & 0.0754 & 0.0811 & 0.0816 & 1.006 & 0.333 & 0.186\\
5   & 0.0600 & 0.0642 & 0.0638 & 0.0686 & 1.070 & 0.667 & 0.442\\
10  & 0.0450 & 0.0537 & 0.0471 & 0.0567 & 1.193 & 0.778 & 0.568\\
30  & 0.0225 & 0.0351 & 0.0230 & 0.0364 & 1.560 & 0.852 & 0.692\\
100 & 0.0082 & 0.0172 & 0.0083 & 0.0175 & 2.097 & 0.878 & 0.759\\
\bottomrule
\end{tabular}
\end{table}

The two edges of the band have a direct meaning. $D_T$ is the budget needed to \emph{keep} a compliant society compliant, and $D_{SN}$ the budget needed to \emph{rescue} a collapsed one with a constant policy; their ratio measures the cost of hysteresis. Hysteresis does not require sanctions that make cheating unprofitable: at the lower edge of the band the expected yearly sanction of a cheater in a compliant society is $FD_T=F(\lambda-\mu)/(1+\lambda F)<1-\mu/\lambda$, smaller than the gain for every $F$, and at $F=F_c$ it equals $\mu/\lambda$, a tenth of the gain for the NHP rates. What matters is not the size of the sanction but how fast its effect fades as compliance falls.

\begin{prop}\label{prop:ratio}
The ratio $\mathcal R(F)=D_{SN}/D_T$ equals $1$ at $F=F_c$ and tends to
\begin{equation}\label{eq:Rinf}
R_\infty=\frac{\lambda^2}{4\mu(\lambda-\mu)}
\end{equation}
as $F\to\infty$. For $\lambda=0.1$, $\mu=0.01$, $R_\infty=2.78$.
\end{prop}

\begin{pf}
For large $F$, $\sqrt{1+F\mu}-\sqrt{F\mu}=1/(\sqrt{1+F\mu}+\sqrt{F\mu})\sim1/(2\sqrt{F\mu})$, so $D_{SN}\sim\lambda/(4F\mu)$, while $D_T\sim(\lambda-\mu)/(\lambda F)$.
\end{pf}

Sanctions lower both budgets, roughly as $1/F$, but they widen the band in relative terms: with a large sanction, rescuing a collapsed society costs almost three times as much as keeping a compliant one, independently of the size of the sanction. The bound depends only on the ratio of imitation to drop-out.

\section{Spontaneous cheating: the cusp in closed form}\label{sec:tau0}

For $\tau_0>0$ the cheater-free state is no longer an equilibrium and the transcritical bifurcation unfolds into a pair of saddle-node bifurcations meeting at a cusp~\cite{Kuznetsov2004}. For $E=0$, $G$ is a cubic polynomial in $y$,
\[
G(y)=-Ay^3+(2A-\lambda+D)\,y^2+(\lambda-D-\mu-\tau_0-A)\,y+\tau_0,\qquad A:=\lambda FD,
\]
and the whole bifurcation set is explicit.

\begin{prop}\label{prop:cusp}
Let $E=0$ and $\tau_0>0$. The cusp point is
\begin{equation}\label{eq:cusp}
y_c=\frac{\tau_0^{1/3}}{\tau_0^{1/3}+\mu^{1/3}},\qquad
A_c=\bigl(\tau_0^{1/3}+\mu^{1/3}\bigr)^3,\qquad
D_c=\lambda-A_c(2-3y_c),\qquad F_{\rm cusp}=\frac{A_c}{\lambda D_c},
\end{equation}
and the fold curves are given parametrically, in terms of the cheater fraction $y\in(0,1)$ at the fold, by
\begin{equation}\label{eq:foldparam}
D(y)=\lambda-\frac{2\mu}{1-y}-\frac{\tau_0\,(1-2y)}{y^2},\qquad
F(y)=\frac{1}{\lambda D(y)}\Bigl(\frac{\mu}{(1-y)^2}+\frac{\tau_0}{y^2}\Bigr).
\end{equation}
As $\tau_0\to0$, $(y_c,D_c,F_{\rm cusp})\to(0,\lambda-2\mu,F_c)$.
\end{prop}

\begin{pf}
At a cusp the cubic has a triple root, $G(y)=-A(y-y_c)^3$. Matching coefficients gives $2A-\lambda+D=3Ay_c$, $\lambda-D-\mu-\tau_0-A=-3Ay_c^2$ and $\tau_0=Ay_c^3$. Eliminating $\lambda-D$ between the first two, $A(1-3y_c+3y_c^2)=\mu+\tau_0$; since $1-3y+3y^2=(1-y)^3+y^3$, dividing by the third gives $\bigl((1-y_c)/y_c\bigr)^3=\mu/\tau_0$, whence (\ref{eq:cusp}). For the folds, write $G=\tau_0c+y\,(-Ac^2+uc-\mu)$ with $u=\lambda-D$, $c=1-y$. The conditions $G=0$ and $G'(y)=0$ are linear in $(A,u)$; subtracting $G/y$ from $G'$ gives $y(2Ac-u)=\tau_0/y$, and then $G=0$ gives $Ac^2=\mu+\tau_0c^2/y^2$, which is (\ref{eq:foldparam}).
\end{pf}

Spontaneous cheating shifts the band to larger budgets and raises the sanction needed for bistability: $F_{\rm cusp}=2.44$, $4.68$ and $13.3$ for $\tau_0=10^{-4}$, $10^{-3}$ and $10^{-2}$, against $F_c=1.25$ for $\tau_0=0$. For $F=30$ the band is $\vphi\in[0.0274,0.0372]$ for $\tau_0=10^{-3}$ and $[0.0399,0.0445]$ for $\tau_0=10^{-2}$, against $[0.0230,0.0364]$ for $\tau_0=0$ (Fig.~\ref{fig:phase}).

\section{How long does a rescue take?}\label{sec:rescue}

\begin{figure}[t]
\centering
\includegraphics[width=0.7\textwidth]{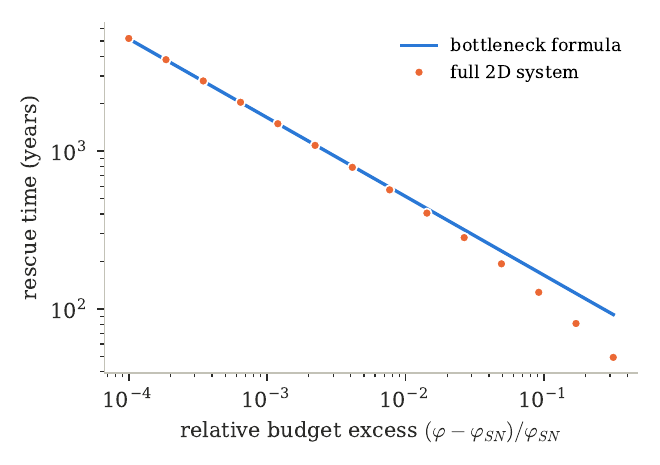}
\caption{Rescue time from the fold state to half its cheater fraction, versus the relative excess of the budget share over $\vphi_{SN}$ ($F=30$, $s=0$, $\tau_0=0$). Line: bottleneck formula (\ref{eq:rescue}); dots: integration of the full system (\ref{eq:model}).}
\label{fig:rescue}
\end{figure}

Just above the fold, the collapsed state has disappeared but the dynamics still slows down where it was: the trajectory has to pass a bottleneck, the ``ghost'' of the saddle-node~\cite{Strogatz2015}. Let $\tau_0=0$, $E=0$, and suppose the society sits at the fold state $y_v=1-c_v$ when the policy is raised from $D_{SN}$ to $D_{SN}+\Delta D$. Expanding (\ref{eq:G}) to leading order around $(y_v,D_{SN})$, using $h(c_v)=h'(c_v)=0$ and $h''=-2\lambda FD$,
\[
\dot y\approx-y_v\Bigl[c_v(1+\lambda Fc_v)\,\Delta D+\lambda FD_{SN}\,(y-y_v)^2\Bigr].
\]
Integrating from $y_v$ gives the time to leave the bottleneck. Along the slow manifold $w=c/\rho+\dot y/\rho^2+\dots$, the lag of the wealth slows the cheater dynamics by the factor $1+\kappa$, $\kappa=-J_{12}/\rho^2$, with $J_{12}$ as in the proof of Proposition~\ref{prop:reduction}. The result is
\begin{equation}\label{eq:rescue}
T\simeq\frac{\pi\,(1+\kappa)}{2\,y_v\sqrt{\lambda FD_{SN}\,c_v(1+\lambda Fc_v)\,\Delta D}},\qquad
\kappa=\frac{k(1-s)\vphi_{SN}\,y_v\,(1+\lambda Fc_v)}{(\theta+\vphi_{SN})^2}.
\end{equation}
Figure~\ref{fig:rescue} compares (\ref{eq:rescue}) with the time taken by the full system (\ref{eq:model}) to reach $y_v/2$, for $F=30$ ($y_v=0.692$, $\kappa=0.045$): they agree to within $1\%$ for excess budgets below $0.3\%$ and to within $4\%$ below $1\%$; further away the quadratic approximation overestimates the time. The rescue time diverges as $\Delta D^{-1/2}$. With a budget $1\%$ above $\vphi_{SN}$ the rescue takes about 490 years, and with $10\%$ about 120 years. Since the rates of the model are per year, this shows that a constant policy only slightly above the threshold is no rescue at all on human time scales, and that faster rescues require policies that vary in time.

\section{Optimal intervention design}\label{sec:optimal}

The closed forms of the previous sections turn the bifurcation analysis into a set of design problems. We consider two: how to split a budget between enforcement and education so as to keep a society compliant at minimum cost, and how to choose the sanction when both cost and resilience matter. Throughout this section $\tau_0=0$, so that maintaining compliance means keeping the cheater-free state stable.

\subsection{Minimum maintenance budget}\label{sec:split}

Write $x=\vphi/(\theta+\vphi)\in[0,1)$ for the fraction of total public spending devoted to fighting cheating, and $s\in[0,1]$ for the fraction of it allocated to education, the rest going to enforcement; then $D=k(1-s)x$ and $E=\eta sx$. In the two-class model education acts only through the moral cost of cheating, measured by $\eta$. By Theorem~\ref{thm:criterion} and (\ref{eq:criterion2}), the cheater-free state is asymptotically stable if and only if
\begin{equation}\label{eq:H}
H(s,x):=\frac{\lambda}{1+\eta sx}-\mu-k(1-s)x(1+\lambda F)<0 .
\end{equation}
The decision maker looks for the smallest effort that keeps compliance,
\begin{equation}\label{eq:design1}
x_T=\inf\bigl\{x\in[0,1):\ H(s,x)<0\ \text{for some } s\in[0,1]\bigr\},
\end{equation}
the boundary $H=0$ being the threshold at which compliance is lost.

\begin{prop}\label{prop:split}
Let $\tau_0=0$, $\rho=\lambda/\mu$ and
\begin{equation}\label{eq:Q}
Q:=\sqrt{\frac{\lambda\eta}{k(1+\lambda F)}} .
\end{equation}
The solution of (\ref{eq:design1}) is:
\begin{enumerate}
\item[(i)] if $Q\le1$ (enforcement only): $s^*=0$ and $x_T=(\lambda-\mu)/(k(1+\lambda F))$;
\item[(ii)] if $1<Q\le\rho$ (mixed programme):
\begin{equation}\label{eq:split}
s^*=\frac{Q-1}{2Q-1-Q^2/\rho},\qquad
x_T=\frac{\lambda(2Q-1)/Q^2-\mu}{k(1+\lambda F)} ;
\end{equation}
\item[(iii)] if $Q>\rho$ (education only): $s^*=1$ and $x_T=(\lambda-\mu)/(\mu\eta)$.
\end{enumerate}
\end{prop}

\begin{pf}
$H$ decreases in $x$ for every $s$, so $x_T$ is the value of $x$ at which $\min_s H=0$. $H$ is convex in $s$ and $\partial_sH=x\bigl(k(1+\lambda F)-\lambda\eta/(1+\eta sx)^2\bigr)$, which is positive for all $s\ge0$ if $Q\le1$, giving (i), and otherwise vanishes at $1+\eta sx=Q$. Substituting $1+\eta sx=Q$ and $k(1+\lambda F)/\eta=\lambda/Q^2$ into $H=0$ gives $x_T$ in (\ref{eq:split}) and $\eta x_T=2Q-1-Q^2/\rho$, hence $s^*=(Q-1)/(\eta x_T)$ as in (\ref{eq:split}); this point lies in $(0,1]$ if and only if $Q\le\rho$. If $Q>\rho$, $H$ decreases in $s$ on $[0,1]$ at the threshold, the minimum is at $s=1$, and $H(1,x)=\lambda/(1+\eta x)-\mu$ vanishes at $x=(\lambda-\mu)/(\mu\eta)$.
\end{pf}

\begin{figure}[t]
\centering
\includegraphics[width=\textwidth]{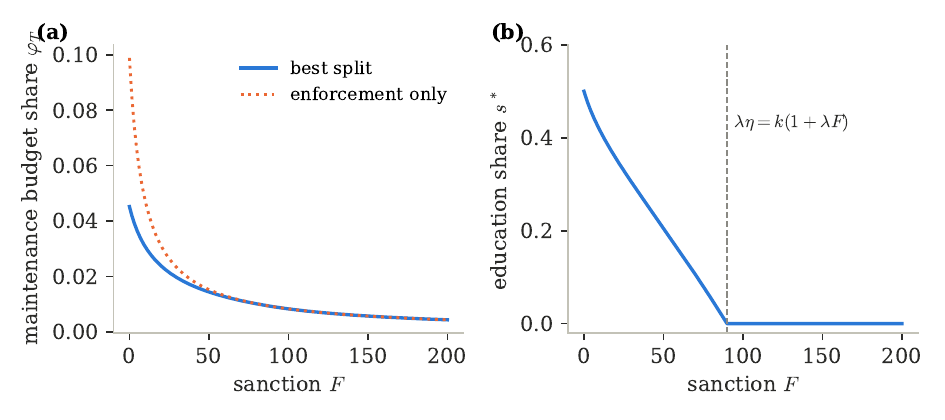}
\caption{(a) Smallest budget share that keeps a compliant society compliant, from (\ref{eq:split}), with the best split (solid) and with enforcement only (dotted). (b) Education share $s^*$ of the best split. $\eta=100$, $\tau_0=0$. The dashed line marks $Q=1$, i.e.\ $\lambda\eta=k(1+\lambda F)$, above which education leaves the cheapest maintenance policy.}
\label{fig:split}
\end{figure}

The parameter $Q^2$ compares the marginal effect of education on the stability of compliance, $\lambda\eta$ at $s=0$, with the marginal effect of enforcement, $k(1+\lambda F)$. The optimal policy therefore has three regimes. If $Q\le1$, equivalently $\lambda\eta\le k(1+\lambda F)$, the whole budget should go to enforcement. If $1<Q\le\lambda/\mu$, the optimal programme is mixed, with the education share (\ref{eq:split}); since $\lambda(2Q-1)/Q^2<\lambda$ for $Q\neq1$, education then always lowers the maintenance budget. If $Q>\lambda/\mu$, education is so effective that the whole budget should go to it; the education share (\ref{eq:split}) reaches $1$ continuously at $Q=\lambda/\mu$. The sanction multiplies the marginal value of enforcement by $1+\lambda F$, because each detection also changes the incentives of the imitators, and so it crowds education out: for $\eta=100$, education is part of the cheapest policy only for $F<(\lambda\eta/k-1)/\lambda=90$ (Fig.~\ref{fig:split}). Without sanctions and with $\eta=100$, the cheapest maintenance policy uses $s^*=0.5$ and $\vphi=0.0452$, less than half the budget needed with enforcement alone ($\vphi=0.099$).

At the optimal split, $\sigma<0$ in (\ref{eq:criterion}) for every $F$ in Fig.~\ref{fig:split}: the cheapest way to keep a society compliant places it inside a bistable band, so that a society kept at this budget does not recover from a shock that pushes it into the cheating state. Cost minimization alone thus leads to fragile policies, which motivates the second design problem.

\subsection{Sanctions: efficiency versus resilience}\label{sec:resilience}

Both edges of the bistable band, $D_T$ and $D_{SN}$ in (\ref{eq:band}), decrease with $F$, so stronger sanctions lower both the budget that keeps a compliant society compliant and the budget that rescues a collapsed one. With $F=30$, for instance, the rescue budget ($\vphi_{SN}=0.036$) is less than half the maintenance budget at $F=F_c$ ($\vphi_T=0.087$). In budgetary terms, sanctions are efficient. They are not free of risk, however: for $F\le F_c$ compliance is lost and regained at the same budget, whereas for $F>F_c$ a shock can push a compliant society into the cheating state, from which only the larger budget $D_{SN}$ rescues it, and slowly (Section~\ref{sec:rescue}).

\begin{figure}[t]
\centering
\includegraphics[width=0.7\textwidth]{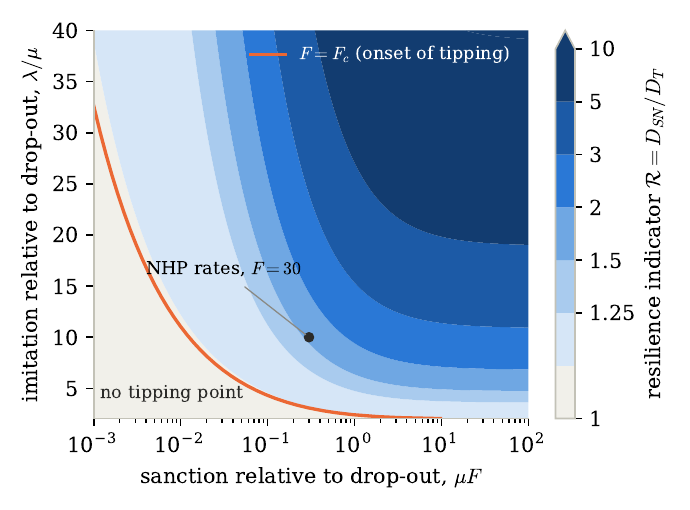}
\caption{Resilience indicator (\ref{eq:resilience}) as a function of the sanction relative to drop-out, $\mu F$, and of the imitation rate relative to drop-out, $\lambda/\mu$ ($\tau_0=0$, $E=0$). Left of the orange line ($F<F_c$) there is no tipping point and $\mathcal R=1$. The dot marks the NHP rates with $F=30$.}
\label{fig:resilience}
\end{figure}

We measure this fragility with the \emph{resilience indicator}
\begin{equation}\label{eq:resilience}
\mathcal R:=\frac{D_{SN}}{D_T}=\frac{(1+\rho m)\bigl(\sqrt{1+m}-\sqrt m\bigr)^2}{1-1/\rho}\quad (m>m_c),\qquad \mathcal R=1\quad (m\le m_c),
\end{equation}
which depends only on two dimensionless groups, the imitation rate relative to drop-out, $\rho=\lambda/\mu$, and the sanction relative to drop-out, $m=\mu F$; the threshold $m_c=\mu F_c=1/[\rho(\rho-2)]$ is the onset of tipping. For small budgets $\mathcal R$ is also the ratio of the budget shares $\vphi_{SN}/\vphi_T$. A value $\mathcal R=1$ means a resilient system without a tipping point; $\mathcal R\gg1$ means a fragile system whose recovery is expensive. For every $\rho$ in Fig.~\ref{fig:resilience}, $\mathcal R$ increases with $m$ towards $R_\infty=\rho^2/[4(\rho-1)]$ (Proposition~\ref{prop:ratio}), so fragility grows with the sanction and, above all, with the speed of imitation relative to drop-out. For the NHP rates ($\rho=10$) and $F=30$, $\mathcal R=1.56$.

The trade-off can be posed as a design problem. With enforcement only, a decision maker who accepts a recovery budget at most $\mathcal R_{\max}$ times the maintenance budget, and who cannot impose sanctions above a legal cap $F_{\max}$, solves
\begin{equation}\label{eq:design2}
\min_{F}\ D_T(F)\quad\text{subject to}\quad \mathcal R(F)\le\mathcal R_{\max},\quad 0\le F\le F_{\max}.
\end{equation}
Since $D_T$ decreases and $\mathcal R$ increases with $F$, the solution is the largest admissible sanction, $F^*=\min\{F_{\max},\,\mathcal R^{-1}(\mathcal R_{\max})\}$, where $\mathcal R^{-1}$ is the inverse of (\ref{eq:resilience}) on $F>F_c$; a decision maker who wants no tipping point at all ($\mathcal R_{\max}=1$) chooses $F^*=\min\{F_{\max},F_c\}$. Table~\ref{tab:design} gives the solution for the NHP rates. Requiring full resilience costs $\vphi_T=0.087$; tolerating a recovery budget twice the maintenance budget allows $F^*=80$ and cuts the maintenance budget to $\vphi_T=0.010$. Where on this curve to stand depends on how likely shocks to compliance are and on how costly a prolonged period of widespread cheating would be, information that lies outside the model but that (\ref{eq:resilience}) and (\ref{eq:design2}) allow a decision maker to weigh explicitly.

\begin{table}[t]
\centering
\caption{Solution of the resilience-constrained design problem (\ref{eq:design2}) without legal cap ($F_{\max}=\infty$), enforcement only, $\tau_0=0$, $\lambda=0.1$, $\mu=0.01$.}
\label{tab:design}
\begin{tabular}{ccccc}
\toprule
$\mathcal R_{\max}$ & $F^*$ & $\vphi_T$ (maintain) & $\vphi_{SN}$ (rescue)\\
\midrule
1 (no tipping point) & 1.25 & 0.0870 & 0.0870\\
1.10 & 6.2 & 0.0587 & 0.0650\\
1.25 & 12.5 & 0.0417 & 0.0526\\
1.50 & 25.9 & 0.0257 & 0.0390\\
2.00 & 80.0 & 0.0101 & 0.0204\\
\bottomrule
\end{tabular}
\end{table}

\section{Sensitivity and robustness}\label{sec:verification}

The parameters of the model, in particular the imitation and drop-out rates, are hard to estimate, so it matters how strongly the budgets depend on them. We also check the two approximations behind the analysis: the limit of a large population, used to derive (\ref{eq:model}) in Section~\ref{sec:derivation}, and the description of the policy as constant in time.

\subsection{Sensitivity of the thresholds}

Let $S^{X}_{p}=(p/X)\,\partial X/\partial p$ denote the elasticity of a quantity $X$ with respect to a parameter $p$: a $1\%$ change in $p$ changes $X$ by $S^X_p\,\%$. With enforcement only and $\tau_0=0$, (\ref{eq:band}) gives
\begin{equation}\label{eq:elasticities}
\begin{aligned}
S^{D_T}_\lambda&=\frac{\lambda}{\lambda-\mu}-\frac{\lambda F}{1+\lambda F},\qquad S^{D_T}_\mu=-\frac{\mu}{\lambda-\mu},\qquad S^{D_T}_F=-\frac{\lambda F}{1+\lambda F},\\
S^{D_{SN}}_\lambda&=1,\qquad S^{D_{SN}}_\mu=S^{D_{SN}}_F=-\sqrt{\frac{\mu F}{1+\mu F}},
\end{aligned}
\end{equation}
and $S^{\mathcal R}_p=S^{D_{SN}}_p-S^{D_T}_p$; the elasticities of the budget shares $\vphi$ follow by multiplying by $k(1-s)/(k(1-s)-D)$, which is close to one for small budgets. With the education channel, (\ref{eq:split}) gives $S^{x_T}_\eta=\lambda(1-Q)/[Q^2k(1+\lambda F)x_T]$ in the mixed regime $1<Q\le\lambda/\mu$. Table~\ref{tab:sensitivity} lists the values for the NHP rates.

\begin{table}[t]
\centering
\caption{Elasticities (\ref{eq:elasticities}) of the maintenance budget $D_T$, the rescue budget $D_{SN}$ and the resilience indicator $\mathcal R$ for $\lambda=0.1$, $\mu=0.01$, $F=30$ (enforcement only, $\tau_0=0$). Last row: elasticity of the optimal maintenance effort $x_T$ with respect to $\eta$ ($\eta=100$) for $F=0,10,30$.}
\label{tab:sensitivity}
\begin{tabular}{lccc}
\toprule
& imitation rate $\lambda$ & drop-out rate $\mu$ & sanction $F$\\
\midrule
maintenance budget $D_T$ & $0.36$ & $-0.11$ & $-0.75$\\
rescue budget $D_{SN}$ & $1.00$ & $-0.48$ & $-0.48$\\
resilience indicator $\mathcal R$ & $0.64$ & $-0.37$ & $0.27$\\
\midrule
\multicolumn{4}{l}{optimal maintenance effort $x_T$ w.r.t.\ $\eta$: $-0.50$ ($F=0$), $-0.42$ ($F=10$), $-0.30$ ($F=30$)}\\
\bottomrule
\end{tabular}
\end{table}

The table shows three things. The cost of prevention is controlled mainly by the sanction: the maintenance budget is most sensitive to $F$ and only weakly sensitive to the behavioural rates, so it can be computed reliably even when $\lambda$ and $\mu$ are uncertain. The cost of recovery is controlled mainly by the imitation rate, with elasticity exactly one: a society in which behaviour spreads $10\%$ faster needs a $10\%$ larger rescue budget, so recovery costs are the less certain quantity. And the resilience indicator grows with the imitation rate and with the sanction and decreases with drop-out, in agreement with Fig.~\ref{fig:resilience}. The boundaries of the bistable region are also robust to spontaneous cheating: for $\tau_0=10^{-3}$ the band for $F=30$ shifts from $[0.0230,0.0364]$ to $[0.0274,0.0372]$ (Section~\ref{sec:tau0}).

\subsection{Individual-based simulations}
\begin{figure}[t]
\centering
\includegraphics[width=0.7\textwidth]{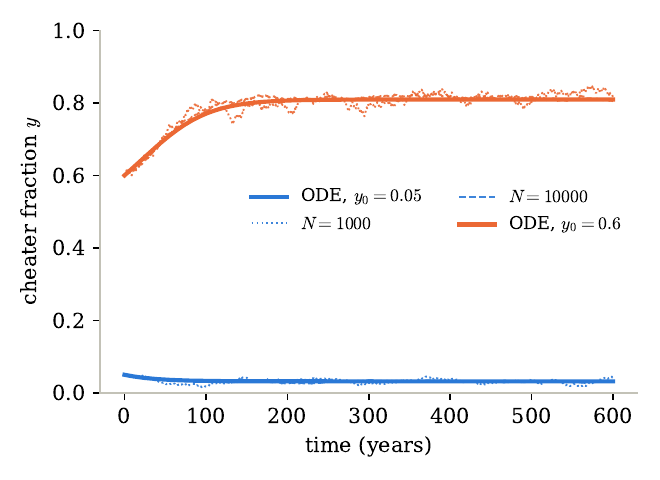}
\caption{Verification of the large-population limit. Thick lines: solution of (\ref{eq:model}) from $5\%$ and $60\%$ of cheaters; thin lines: individual-based simulations with $N=1000$ (dotted) and $N=10\,000$ (dashed). $F=30$, $s=0$, $\tau_0=10^{-3}$, $\vphi=0.032$.}
\label{fig:verification}
\end{figure}

We simulated the individual process of Steps~1 and~2 of Section~\ref{sec:derivation} directly, without the mean-field limit, by $\tau$-leaping with step $\Delta t=0.02$ years: in each step the numbers of payers who start cheating and of cheaters who stop are Poisson variables with the rates of (\ref{eq:expected}), and the wealth is updated with the current number of payers. We used $F=30$, $s=0$, $\tau_0=10^{-3}$ and $\vphi=0.032$, inside the bistable band $[0.0274,0.0372]$ of Section~\ref{sec:tau0}. Figure~\ref{fig:verification} shows that simulations of $N=1000$ and $N=10\,000$ individuals follow the solution of (\ref{eq:model}) and, starting from $5\%$ and $60\%$ of cheaters, settle in the compliant and in the cheating state respectively: the coexistence of the two states is not an artefact of the large-population limit. Over 20 simulations per size, the cheater fraction after 300 years differs on average from the deterministic value ($0.0319$) by less than $10^{-3}$, and its standard deviation decreases from $6.3\times10^{-3}$ for $N=500$ to $1.0\times10^{-3}$ for $N=32\,000$, as expected from fluctuations of order $N^{-1/2}$.

\subsection{A slow budget cycle}

\begin{figure}[t]
\centering
\includegraphics[width=0.7\textwidth]{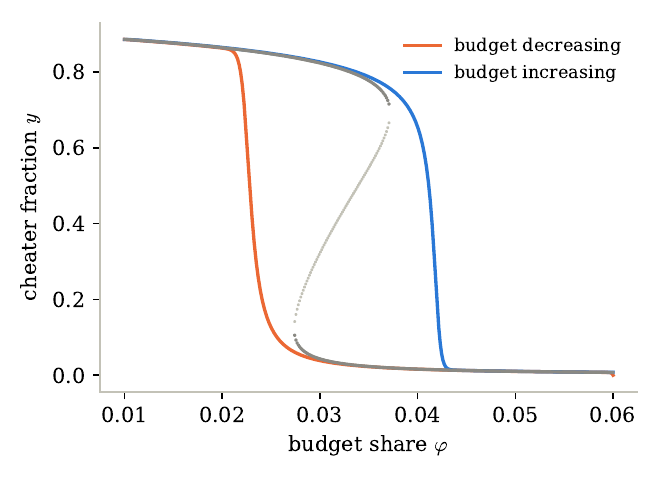}
\caption{Hysteresis under a slow budget cycle: the budget share decreases from $0.06$ to $0.01$ (orange) and increases back (blue) over 3000 years each way. Grey dots: stable (dark) and unstable (light) equilibria for constant budgets. $F=30$, $s=0$, $\tau_0=10^{-3}$.}
\label{fig:cycle}
\end{figure}

To show what hysteresis means for a decision maker, we let the budget share decrease slowly from $\vphi=0.06$ to $0.01$ over 3000 years and then increase back at the same rate ($F=30$, $s=0$, $\tau_0=10^{-3}$). The society collapses when the budget falls to $\vphi\approx0.023$ and recovers only when it rises to $\vphi\approx0.042$ (Fig.~\ref{fig:cycle}). Restoring the budget at which compliance was lost is therefore not enough. Both switches occur beyond the static edges of the band, $0.0274$ and $0.0372$, because the dynamics slows down near the folds (Section~\ref{sec:rescue}): when the budget changes, the effect appears with a delay, so that a collapse becomes visible only after the threshold has already been crossed.

\section{Managerial insights: a decision table}\label{sec:decision}

Table~\ref{tab:decision} collects, for three levels of the sanction, the quantities a decision maker needs: the budget share that keeps a compliant society compliant ($\vphi_{\rm maintain}$, the lower edge of the band), the budget share above which a collapsed society recovers ($\vphi_{\rm rescue}$, the upper edge), the fraction of cheaters in the collapsed state, and the time needed to bring the fraction of cheaters below $5\%$ with a budget $10\%$, $25\%$ and $50\%$ above $\vphi_{\rm rescue}$. All values include spontaneous cheating ($\tau_0=10^{-3}$) and use enforcement only ($s=0$).

\begin{table}[t]
\centering
\caption{Decision table ($\tau_0=10^{-3}$, $s=0$, other parameters as in Section~\ref{sec:model}). $T_{+x\%}$: years needed to bring the fraction of cheaters from the collapsed state below $5\%$ with a budget share $x\%$ above $\vphi_{\rm rescue}$.}
\label{tab:decision}
\begin{tabular}{rcccccccc}
\toprule
$F$ & $\vphi_{\rm maintain}$ & $\vphi_{\rm rescue}$ & $\vphi_{\rm rescue}/\vphi_{\rm maintain}$ & collapsed $y$ & $T_{+10\%}$ & $T_{+25\%}$ & $T_{+50\%}$\\
\midrule
10  & 0.0541 & 0.0581 & 1.07 & 0.72 & 363 & 168 & 90\\
30  & 0.0274 & 0.0372 & 1.36 & 0.84 & 274 & 147 & 87\\
100 & 0.0100 & 0.0179 & 1.79 & 0.87 & 236 & 131 & 80\\
\bottomrule
\end{tabular}
\end{table}

Three practical conclusions follow. First, stronger sanctions lower both budgets markedly, but they widen the gap between them: with $F=100$ a collapsed society needs almost $1.8$ times the budget that would have kept it compliant, approaching the bound $R_\infty=2.78$ of Proposition~\ref{prop:ratio}. Second, the gap is only part of the cost of a collapse. Even with a budget $50\%$ above the rescue threshold, the recovery takes about 80--90 years with the rates of the original model, and a budget just $10\%$ above it takes centuries, so the time lost during a collapse is the dominant cost. Since all rates scale with the imitation rate $\lambda$, these times would be shorter in a society where behaviour spreads faster, but their ratios would not change. Third, because the maintenance budget can be computed in closed form (Proposition~\ref{prop:split}), an agency can monitor how far its current budget is from the threshold and act before compliance is lost, which is much cheaper than acting afterwards.

\section{Discussion}\label{sec:discussion}

The two-class model is small enough to answer exactly a question that richer models address only numerically: what makes the fight against cheaters history-dependent? In this model the answer is the conjunction of two ingredients. The first is a fiscal feedback, which makes enforcement proportional to compliance. The second is a deterrent that acts on incentives, so that cheaters respond to it by changing their behaviour and not only by being removed. Police financed by the tax base, without sanctions, affects cheaters in exactly the same way as imitation itself, and only renormalizes the imitation rate; a budget set as an absolute amount has no feedback at all. In both cases the steady state is unique, which explains the uniqueness that NHP found for their two-class model. Sanctions and education that raises the moral cost of cheating both fade when compliance falls, and both create a tipping point once they exceed the explicit thresholds (\ref{eq:Fc}) and Corollary~\ref{cor:education}. The general criterion (\ref{eq:criterion}) states the condition for any such channel.

The closed forms also give three quantitative messages. The sanction needed for hysteresis is small when cheating spreads mainly by imitation ($F_c=1.25$ years of tax for the NHP rates) and grows with spontaneous cheating (Proposition~\ref{prop:cusp}). The cost of rescuing a collapsed society, relative to keeping a compliant one, is bounded by $\lambda^2/[4\mu(\lambda-\mu)]$, whatever the sanction. And a constant policy just above the rescue threshold is slow, with a time that diverges as the inverse square root of the excess budget.

\paragraph{Implications for decision makers.} Within the limits of the model, the results suggest five practical points for anyone managing a compliance system funded by those who comply. First, the risk of a tipping point is created by linking enforcement to current revenue: an enforcement budget protected from falls in revenue, fixed in absolute terms, cannot be eroded by cheating and removes the hysteresis (Corollary~\ref{cor:money}). Second, the relevant quantity for an enforcement agency is not the level of deterrence but how quickly it weakens when compliance falls, which is what the criterion (\ref{eq:criterion}) measures; fines and education that are effective in a compliant society can make it fragile, because their effect depends on resources that cheating destroys. Third, prevention is cheaper than cure: keeping a society compliant costs a fraction of what it costs to rescue it once it has collapsed, and a rescue with a budget barely above the threshold is too slow to be useful, so a timely response to a drop in compliance matters more than its eventual size. Fourth, the best use of a given budget depends on the sanctions in force: when fines are mild, part of the budget should go to education; when they are strong, enforcement is the more efficient instrument. Fifth, raising sanctions buys efficiency at the price of resilience: above $F_c$ every budget is lower but the system acquires a tipping point, and the resilience indicator (\ref{eq:resilience}) and the design problem (\ref{eq:design2}) make the trade-off quantitative.

\paragraph{Other compliance systems.} Nothing in the analysis is specific to taxes. The model applies whenever compliance is a binary choice that spreads by imitation, deterrence has a mechanical component (removal) and an incentive component (sanctions, moral cost), and the enforcement budget is a share of what the compliant contribute. Fare evasion in public transport, where inspections are paid from ticket revenue, the payment of fees that fund the monitoring of fishing quotas, waste dumping or other environmental regulations, and corruption control financed from the public budget that corruption depletes share this structure. In each case Theorem~\ref{thm:criterion} identifies the instruments that can create a tipping point, and the closed forms give the corresponding maintenance and rescue budgets, once the rates are estimated for the system at hand.

\paragraph{Limitations.} The model has no social structure: education can act only on incentives, and effects of policy on the class composition of the population, present in the three-class model of NHP, are absent (Remark~\ref{rem:threeclass}). Fines are not returned to the common wealth, which would add a further feedback, and the sanction is a fixed multiple of the yearly tax. Policies are constant in time, and the population is well mixed; time-dependent policies and spatial structure are natural extensions. Finally, Theorem~\ref{thm:criterion} is a local statement at full compliance, where only the derivative of the switching rate at zero payoff difference matters; for a nonlinear imitation rule we therefore expect the criterion to hold with $\lambda$ replaced by that derivative, while the closed forms of Sections~\ref{sec:sanctions}--\ref{sec:optimal} use the linear rule. The parameters are those of the original NHP model rather than estimates from tax-compliance data; calibrating the imitation and drop-out rates with longitudinal data on evasion would be the next step towards quantitative predictions.

\section{Conclusions}\label{sec:conclusions}

We have developed an exactly solvable model of compliance in which the resources used against cheating come from those who comply. Deriving the population equations from individual imitation separates the two effects of a deterrent, removal and incentives, and this separation yields a sharp criterion: tipping points require that an instrument acting on incentives be funded by the tax base, while fiscally funded enforcement alone, or any budget fixed in absolute terms, never creates them. For sanctions, all the relevant quantities are explicit: the bistable band, the critical sanction, the cusp, the maintenance and rescue budgets, their ratio, and the slow recovery near the threshold. These closed forms support two design problems with explicit solutions, the cheapest split of a maintenance budget between enforcement and education and the choice of the sanction under a resilience constraint, and a sensitivity analysis shows that prevention costs are controlled mainly by sanctions and recovery costs mainly by the speed of imitation.

Beyond the application to tax compliance, the model provides a general analytical framework for quantifying prevention costs, recovery costs and intervention effectiveness in socio-economic systems subject to fiscal feedback. Because the thresholds and recovery requirements are available in closed form, the framework is particularly suited to the design and comparison of resource-constrained intervention policies, and to the calibration with data on compliance that would turn its qualitative guidelines into quantitative prescriptions.

\appendix
\section{Numerical methods}\label{app:numerics}

Equilibria were located by bracketing the zeros of $G$ on a grid of $3001$--$4001$ points of $[0,1]$ and refining them with Brent's method; stability was obtained from the eigenvalues of the $2\times2$ Jacobian of (\ref{eq:model}) computed by central differences. Bistable bands were obtained by scanning $\vphi$ (600--1200 values) and recording where two stable equilibria coexist; for $\tau_0=10^{-10}$ they agree with (\ref{eq:band}) to within the scan resolution ($\sim2\times10^{-5}$ in $\vphi$), and for $\tau_0=10^{-3}$ with the parametric folds (\ref{eq:foldparam}). Rescue times were computed by integrating (\ref{eq:model}) with LSODA (relative tolerance $10^{-10}$) from the fold state $(y_v,(1-y_v)/(\theta+\vphi_{SN}))$ until $y=y_v/2$. The maintenance optimum (\ref{eq:split}) was checked against a direct numerical minimization over $s$. The Python code that reproduces all results and figures is available from the corresponding author on request.

\section*{CRediT authorship contribution statement}
\textbf{Luca Meacci:} Conceptualization, Methodology, Formal analysis, Writing -- review \& editing. \textbf{Francisco J. Mu\~noz:} Conceptualization, Methodology, Formal analysis, Software, Visualization, Writing -- original draft. \textbf{Juan Carlos Nu\~no:} Conceptualization, Methodology, Supervision, Writing -- review \& editing.

\section*{Declaration of competing interest}
The authors declare that they have no known competing financial interests or personal relationships that could have appeared to influence the work reported in this paper.

\section*{Declaration of generative AI and AI-assisted technologies in the manuscript preparation process}
During the preparation of this work the authors used Claude (Anthropic) to help implement the questions formulated by the authors and to revise the text, the code and the figures. After using this tool, the authors reviewed and edited the content as needed and take full responsibility for the content of the published article.

\section*{Data availability}
No data were used for the research described in the article. The Python code that reproduces all results and figures is available from the corresponding author on request.

\section*{Acknowledgements}
This research did not receive any specific grant from funding agencies in the public, commercial, or not-for-profit sectors. L.M. acknowledges the support of the Gruppo Nazionale per la Fisica Matematica (GNFM) of the Istituto Nazionale di Alta Matematica (INdAM). We dedicate this work to Mario Primicerio.


\begin{thebibliography}{99}
\bibitem{Becker1968} G.~S. Becker, Crime and punishment: an economic approach, J. Polit. Econ. 76 (1968) 169--217. \url{https://doi.org/10.1086/259394}

\bibitem{AllinghamSandmo1972} M.~G. Allingham, A. Sandmo, Income tax evasion: a theoretical analysis, J. Public Econ. 1 (1972) 323--338. \url{https://doi.org/10.1016/0047-2727(72)90010-2}

\bibitem{Andreoni1998} J. Andreoni, B. Erard, J. Feinstein, Tax compliance, J. Econ. Lit. 36 (1998) 818--860.

\bibitem{PickhardtPrinz2014} M. Pickhardt, A. Prinz, Behavioral dynamics of tax evasion: a survey, J. Econ. Psychol. 40 (2014) 1--19.

\bibitem{HofbauerSigmund1998} J. Hofbauer, K. Sigmund, Evolutionary Games and Population Dynamics, Cambridge University Press, Cambridge, 1998.

\bibitem{Sigmund2010} K. Sigmund, H. De Silva, A. Traulsen, C. Hauert, Social learning promotes institutions for governing the commons, Nature 466 (2010) 861--863. \url{https://doi.org/10.1038/nature09203}

\bibitem{Nuno2008} J.~C. Nu\~no, M.~A. Herrero, M. Primicerio, A triangle model of criminality, Physica A 387 (2008) 2926--2936.

\bibitem{Nuno2010} J.~C. Nu\~no, M.~A. Herrero, M. Primicerio, Fighting cheaters: how and how much to invest, Eur. J. Appl. Math. 21 (2010) 459--478. \url{https://doi.org/10.1017/S0956792510000094}

\bibitem{Nuno2011} J.~C. Nu\~no, M.~A. Herrero, M. Primicerio, A mathematical model of a criminal-prone society, Discrete Contin. Dyn. Syst. Ser. S 4 (2011) 193--207.

\bibitem{Scheffer2009} M. Scheffer, J. Bascompte, W.~A. Brock, V. Brovkin, S.~R. Carpenter, V. Dakos, H. Held, E.~H. van Nes, M. Rietkerk, G. Sugihara, Early-warning signals for critical transitions, Nature 461 (2009) 53--59. \url{https://doi.org/10.1038/nature08227}

\bibitem{Sandholm2010} W.~H. Sandholm, Population Games and Evolutionary Dynamics, MIT Press, Cambridge, MA, 2010.

\bibitem{Schlag1998} K.~H. Schlag, Why imitate, and if so, how? A boundedly rational approach to multi-armed bandits, J. Econ. Theory 78 (1998) 130--156.

\bibitem{Kurtz1970} T.~G. Kurtz, Solutions of ordinary differential equations as limits of pure jump Markov processes, J. Appl. Probab. 7 (1970) 49--58.

\bibitem{TaylorJonker1978} P.~D. Taylor, L.~B. Jonker, Evolutionary stable strategies and game dynamics, Math. Biosci. 40 (1978) 145--156. \url{https://doi.org/10.1016/0025-5564(78)90077-9}

\bibitem{Perko2001} L. Perko, Differential Equations and Dynamical Systems, 3rd ed., Springer, New York, 2001.

\bibitem{Kuznetsov2004} Yu.~A. Kuznetsov, Elements of Applied Bifurcation Theory, 3rd ed., Springer, New York, 2004.

\bibitem{Strogatz2015} S.~H. Strogatz, Nonlinear Dynamics and Chaos, 2nd ed., Westview Press, Boulder, CO, 2015.
\end{thebibliography}
\end{document}